\documentclass[12pt]{article}
\usepackage[cp1251]{inputenc}
\usepackage{amssymb}
\usepackage{amsfonts}
\usepackage{amsmath}
\def\bee{\begin{eqnarray}}
\def\bes{\begin{eqnarray*}}
\def\eee{\end{eqnarray}}
\def\ees{\end{eqnarray*}}

\begin{document}

\noindent {\bf Free Malcev algebra of rank three.}\\
Alexandr Ivanovich Kornev \\  
{Centro de Matem\'atica, Computa\c c\~ao e Cogni\c c\~ao,\\  
Universidade Federal do ABC, Avenida dos Estados, 5001,  \\
Bairro Bangu, Santo Andr\'e, SP, Brasil, CEP 09.210-580\\  
phone +55(11)49968315 }\\
{\it e-mail: alexandr.kornev@ufabc.edu.br}

\title{\bf Free Malcev algebra of rank three.}
\author{A.I.\ Kornev\footnote{Supported by FAPESP grant No. 2008/57680-5.}\\  
Federal University of ABC (UFABC)\\ Brazil}
\maketitle 

\begin{abstract} 
We find a basis for the free Malcev algebra on
three free generators over a field of characteristic zero. The   semiprimity and speciality 
of this algebra are proved. 
Also, we prove the decomposability of this algebra into a subdirect sum of the free Lie algebra of rank
three and the free algebra of rank three of the variety  generated by a simple seven-dimensional
Malcev algebra. 
These results were announced in ~\cite{Kor}.
\\
\end{abstract}
The problem of finding of a basis for a free  algebra is important for different varieties.
For free Malcev algebras this problem is posed by Shirshov in ~\cite[ problem 1.160]{Dn}.
For alternative algebras with three generators similar problem is solved in ~\cite{Il}.
The base of the free Malcev superalgebra on one odd generator was constructed in ~\cite{Sh1odd}.
Recall that a Malcev algebra is called special if it is a subalgebra of a commutator 
algebra $A^{-}$ for some alternative algebra $A.$
The question of the speciality of a Malcev algebras was posed by Malcev in ~\cite{Mal}.
In the present paper we find the basis of the free Malcev algebra
 with three free generators, and prove the speciality of this algebra.
In addition, we prove decomposition of this algebra into a subdirect sum of free Lie algebra of rank
three and  the free algebra of rank three of the variety  generated by a simple seven-dimensional Malcev algebra.\\
Shestakov ~\cite{Sh2} in 1977 proved that a free Malcev algebra of $n>8$
generators over commutative ring $\Phi$ is not semiprime provided $7!\ne 0$ in $\Phi$.
Filippov ~\cite{Fil4} in 1979 then proved that in fact a free Maltsev $\Phi$-algebra of $n>4$ generators 
is not semiprime if $6\Phi \ne 0.$ We prove that the free Malcev algebra of rank three is semiprime.

For brevity, we omit the brackets in  terms of the  
type  $(...((x_1x_2)x_3)...)x_n.$ In addition, the products of the form 
$axx...x$ we denote by $ax^n.$\\
A linear algebra $M$ over a field $F$, which satisfies
following identities 
 $$x^2=0,$$    $$J(x,y,xz)=J(x,y,z)x,$$ where $J(x,y,z)=xyz+zxy+yzx$, is called a Malcev algebra.
 In what follows, the characteristic of $F$ is assumed to be zero. Let $R_a$ be the operator of right multiplication by
 an element $a$ of the algebra $M$, and let $R(M)$ be the algebra generated by all $R_a$.
 We will adhere to  following notations:\\
                           $$L_{a,b}=1/2(R_aR_b+R_bR_a),$$
                           $$[R_a,R_b]=R_aR_b-R_bR_a,$$
                           $$(a,b,c)=(ab)c-a(bc).$$
   In addition, by $\Delta_a^i(b)$  we denote the operator of partial linearization defined in ~\cite[Chapter 1, \S 4]{zsss},
   by $Z(M)$ we denote the Lie  center of Malcev algebra  $M$ and
  by $G(a,b,c,d)$ we denote the function defined in  ~\cite {sag} as follows: $$G(a,b,c,d)=J(ab,c,d)-bJ(a,c,d)-J(b,c,d)a.$$
 Let $X=\{x,y,z\}$ and $M_3$ be the free Malcev algebra on the set of  free generators $X$. 
For brevity, the expressions of the form $G(...(G(G(t,x,y,z),x,y,z)...),x,y,z)$ will be denoted as $tG^n$. \\
  By $Alt[X]$ denote the free alternative algebra generated by the set of free generators $X$ and by $Ass[X]$ 
denote the free associative algebra, generated by the set of free generators $X$. 
Furthermore, for $a,b \in Alt[X]$ we denote $a\circ b=1/2(ba+ab)$ 
   and  for $a\in Alt[X]$ denote by $R^+_a$ the
     operator in the algebra $Alt[X]$ defined as follows:
   $xR^+_a=x\circ a$ for any $x\in Alt[X]$.
   Now $L^+_{a,b}$, is the operator in the algebra  $Alt[X]$, defined as follows:
   $L^+_{a,b}=R^+_aR^+_b-R^+_{a\circ b}$ for any $x\in Alt[X]$ and $a,b \in Alt[X]$. If $B$ 
   is an alternative algebra, then by $B^{-}$ denote the commutator algebra of the algebra $B$.\\
   The main result:\\
   {\bf Theorem.} Let $M_3$ be the free Malcev algebra with free generators $X=\{x,y,z\}$
    Let $\textbf{U}=\{J(x,y,z) G^k L_{x,x}^l
L_{y,y}^m L_{z,z}^n L_{x,y}^p L_{x,z}^q L_{y,z}^r \mid k,l,m,n,p,q,r
\in \mathbb{N}\cup \{0\}\}$. Then the set of the vectors
$\textbf{U}\cup \textbf{U}x\cup \textbf{U}y\cup \textbf{U}z\cup
\textbf{U}xy\cup \textbf{U}xz\cup \textbf{U}yz$ forms a basis of the space $J(M_3,M_3,M_3).$ 
Furthermore, $M_3$ is special.\\

Any Malcev algebra $M$ and the algebra $R(M)$ satisfy the following identities:\\
\begin{equation}\label{a7}
uL_{a,b}t=utL_{a,b}+uL_{at,b}-uL_{a,tb}.
\end{equation}

\begin{equation}\label{a1}
(ab)(cd)=acbd+dacb+bdac+cbda,
\end{equation}

\begin{equation}\label{a4}
G(t,a,b,c)=2(J(ta,b,c)+J(t,a,bc)),
\end{equation}

\begin{equation}\label{a3}
G(t,a,b,c)=2/3(J(t,b,c)a+J(a,t,c)b+J(a,b,t)c-J(a,b,c)t),
\end{equation}

\begin{equation}\label{a2}
3J(wa,b,c)=J(a,b,c)w-J(b,c,w)a-2J(c,w,a)b+2J(w,a,b)c,
\end{equation}

\begin{equation}\label{a5}
J(J(a,b,c),a,b)=-3J(a,b,c)(ab),
\end{equation}

\begin{equation}\label{a6}
J(c,ba^{2k-1},b)=-J(c,b,a)a^{2k-2}b, (k\geq 1),
\end{equation}
The  identity  \eqref{a7} is the identity (12) of ~\cite[\S 1]{fil}
rewritten in our notation.  The identities  \eqref{a1}, \eqref{a4}, \eqref{a3}, \eqref{a2} and \eqref{a5} are proved in ~\cite{sag}   and \eqref{a6} is
the identity (5) of ~\cite[\S 1]{fil}.
Moreover, from the identity \eqref{a3} it is clear that the function $G$ is  skew-symmetric for any two arguments. \\

 {\bf Lemma 1.} Any Malcev algebra $M$ satisfies the following identities:\\

\begin{equation}\label{a10}
 J(a,b,c)L_{b,b}^ka=J(a,b,c)aL_{b,b}^k,
\end{equation}
\begin{equation}\label{a11}
 J(a,b,c)L_{a,a}^kL_{b,b}^l=J(a,b,c)L_{b,b}^lL_{a,a}^k,
\end{equation}
$$(ta)J(a,b,c)=$$
 \begin{equation}\label{a8}
-1/2J(a,t,c)ab+1/2J(a,t,b)ac-J(b,c,ta)a-3/2J(a,t,cb)a,
\end{equation}
\begin{equation}\label{a9}
J(a,b,tac)=-1/2J(a,t,c)[R_a,R_b]+J(a,b,t)L_{a,c}.
\end{equation}

{\sc Proof}.  The identity \eqref{a10} easy follows from the identity (21)  of ~\cite{Kuz}.\\

The identiy \eqref{a11}. We apply  the  identities \eqref{a6} and \eqref{a10}.
$$J(a,b,c)L_{a,a}^{k}L_{b,b}^{l}=-J(ba^{2k+1},b,cb^{2l-1})=-J(a,b,cb^{2l-1})L_{a,a}^{k}b=$$
$$-J(a,b,cb^{2l-1})bL_{a,a}^{k}=J(a,b,c)L_{b,b}^{l}L_{a,a}^{k}.$$

The identiy \eqref{a8}. Applying the operator $\Delta_b^1(h)$ to the identity $bJ(a,b,c)=J(a,b,cb)$ we obtain
$$hJ(a,b,c)=J(a,h,c)b+J(a,h,cb)+J(a,b,ch).$$ From the identity \eqref{a2}:
$$hJ(a,b,c)=1/3J(a,h,c)b+J(a,h,cb)+1/3J(a,b,h)c-2/3J(b,c,h)a+1/3hJ(a,b,c).$$
That is, $$hJ(a,b,c)=1/2J(a,h,c)b-J(b,c,h)a+1/2J(a,b,h)c+3/2J(a,h,cb).$$
Replacing now $h$ by $ta$ we obtain the identity \eqref{a8}.\\

  The identiy \eqref{a9}. Applying twice the first the identity \eqref{a2} to the identity \eqref{a8}, we obtain the identity
$$(ta)J(a,b,c)=-1/2J(a,t,c)ab-1/6J(a,b,t)ca-1/6J(a,t,c)ba-1/6J(t,b,c)aa+$$
 $$+2/3J(a,b,c)ta-1/2J(a,b,t)ac.$$ From the identity \eqref{a2} we obtain
$$J(a,b,tac)=2/3J(ta,c,a)b-2/3J(b,ta,c)a-1/3J(a,b,ta)c+1/3J(c,a,b)(ta)=$$
$$-2/3J(t,c,a)ab+4/9J(a,b,t)ca+4/9J(a,t,c)ba-2/9J(t,b,c)aa+$$
$$2/9J(a,b,c)ta+1/3J(a,b,t)ac-1/3(ta)J(a,b,c).$$
 Applying now the previous identity, we obtain \eqref{a9}.$\blacksquare$ \\


 {\bf Lemma 2.} For an arbitrary polynomial $f$ of degree $n$ from the subalgebra
 generated by the elements $a$ and $b$ of Malcev algebra $M$,  the following equalities are true:\\
1) $J(a,b,J(a,f,c))+(-1)^nJ(a,f,J(a,b,c))=0,$\\
2) $J(a,b,J(f,b,c))+(-1)^nJ(f,b,J(a,b,c))=0,$\\
3) $J(a,b,fc)=(-1)^nfJ(a,b,c).$\\

{\sc Proof}. We shall prove  all three statements by induction on $n$. For $n =1$ all identities are obvious.
We suppose they are correct when $ n = k $, and prove them for $ n = k +1 $. 
Since  $M$ is a binary Lie algebra,  we can assume that $f = f_1a$
or $ f = f_1b $. \\
1) If $f=f_1a$, the proof is obvious. Let $f=f_1b$. Applying
operator $\Delta_c^1(f_1)$ to the identity $J(a,cb,c)=J(a,b,c)c$ we obtain
$$J(a,f_1b,c)=J(a,f_1,cb)+J(a,b,c)f_1.$$ Next, using this equation and the induction hypothesis, obtain
$$J(a,b,J(a,f_1b,c))+(-1)^{k+1}J(a,f_1b,J(a,b,c))= J(a,b,J(a,f_1,cb))+J(a,b,J(a,b,c)f_1)+$$
$$+(-1)^{k+1}J(a,f_1,J(a,b,c)b)+(-1)^{k+1}J(a,b,J(a,b,c))f_1=0$$
2) The proof is similar to  1).\\
3) Let $f=f_1a$.
$$J(a,b,f_1ac)=J(a,b,J(f_1,a,c)-cf_1a-acf_1)=$$
$$-J(a,b,J(a,f_1,c))+(-1)^{k+1}(f_1J(a,b,c)a+J(a,b,c)af_1)=$$
$$-J(a,b,J(a,f_1,c))+(-1)^{k+1}(J(J(a,b,c),a,f_1)-af_1J(a,b,c))=(-1)^{k+1}f_1aJ(a,b,c).$$\\
If $f=f_1b$, then the arguments are similar and use the equality
2).$\blacksquare$\\

{\bf Corollary 1.} Under the conditions of Lemma 2 the following equalities hold :\\
\begin{equation}\label{a12}
 J(a,f,c)b-J(b,f,c)a=\frac {3(-1)^n+1}{2}fJ(a,b,c),
\end{equation}
 \begin{equation}\label{a13}
 J(a,fb,c)=-J(a,f,c)b+(-1)^nfJ(a,b,c).
\end{equation}

{\sc Proof}. From Lemma 2:
$J(a,b,fc)=(-1)^nfJ(a,b,c).$ From the identity \eqref{a2}: 
$$J(a,b,fc)=2/3J(f,c,a)b-2/3J(b,f,c)a-1/3J(a,b,f)c+1/3J(c,a,b)f.$$ 
Combining, we obtain the first equality. The second
equality follows from the first after the application 
of the identity \eqref{a2} to $J(a,fb,c)$.
$\blacksquare$\\

 {\bf Proposition 1.} Let $u\in J(M_3,M_3,M_3)$.
There exist $\alpha_i$ from $F$ such that
$$u=\sum_{i}\alpha_iJ(x,y,z)x_{i1}x_{i2}...x_{ik_i},$$
where  $x_{ij}\in X, k_i\in \mathbb{N}\cup \{0\}.$\\
{\sc Proof.} In any Malcev algebra  $M$  the following identity is fulfilled:
$$J(a,b,c)abv=J(a,b,v)cba+J(a,b,c)vba+J(a,b,v)bca+$$
\begin{equation}
+J(a,b,c)bva-J(a,b,v)abc-J(a,b,c)avb-J(a,b,v)acb.\tag{a}
\end{equation}
Applying the operator $\Delta_b^1(c)$ to identity \eqref{a10} for $k=1$ obtain:
$$J(a,b,c)cba+J(a,b,c)bca=J(a,b,c)acb+J(a,b,c)abc.$$ 
After the application of the operator $\Delta_c^1(v)$ we obtain the desired identity.\\

We prove the proposition  by induction on the degree  $d$ of the
element $u$. Corollary 2 ~\cite[Chapter 1, \S 3]{zsss} implies
the homogeneity of any variety of Malcev algebras over a field  $F$.
Therefore, the proposition is valid for $d=3$.  Assume the proposition  true for all $d\leq k$. 
 From the identity \eqref{a2} and the induction hypothesis it follows that $u$ can be written in the form
$$u=\sum_{i}\alpha_iJ(x,y,z)x_{i1}x_{i2}...x_{ik_i}v_i,$$ where $x_{ij}\in X, k_i\in \mathbb{N}\cup \{0\},v_i\in M_3.$
It is obvious that the elements 
$\overline{v}_i=v_i+J(M_3,M_3,M_3)$ of the Lie algebra  $M_3/J(M_3,M_3,M_3)$ are linear combinations of the 
monomials of the form 
$\overline{y}_{i1} \overline{y}_{i2}...\overline{y}_{il}$, where
$y_{ij}\in X$ and  $l\in \mathbb{N}$. Hence $v_i$
can be written as $v_i=y_{i1}y_{i2}...y_{il}+u_i$, where
$u_i\in J(M_3,M_3,M_3)$. From  homogeneity it follows that the degrees of the elements
$u_i$ does not exceed $k$, therefore, the induction hypothesis
implies that $v_i$ can be represented as a linear combination of the 
monomials of the form $t_{i1}t_{i2}...t_{is}$, where $t_{ij}\in X$ and
$s\leq k$. Thus, to prove the proposition it is sufficient
show that expressions of the form $J(x,y,z)x_{1}x_{2}...x_{r}(wt)$, where
$x_{j},t\in X$ and $w\in M_3$, belong to the subspace
generated by the set $J(M_3,M_3,M_3)X$. We prove this by 
induction on $r$. From the identity \eqref{a8} it follows that this is true for
$r=0$. Let $r=1$. There are two cases.\\
1. $x_1=t$. Without loss of generality we can assume that $x_1=t=x$.
We have
$$J(x,y,z)x(wx)=J(J(x,y,z),x,wx)-(wx)J(x,y,z)x+(wxx)J(x,y,z)=$$
$$-J(J(x,y,z),x,w)x-(wx)J(x,y,z)x+(wxx)J(x,y,z).$$ From the identity \eqref{a8} 
it follows that the third term, and therefore all expressions lie in $J(M_3,M_3,M_3)X$.\\
2. $x_1\ne t$. Without loss of generality, we may assume that $x_1=x, t=y$.
Applying the identities \eqref{a1} and \eqref{a8}, we obtain 
$$J(x,y,z)x(wy)=J(x,y,z)wxy+yJ(x,y,z)wx+xyJ(x,y,z)w+wxyJ(x,y,z)=$$
$$J(x,y,z)wxy+yJ(x,y,z)wx+wxyJ(x,y,z)+1/2J(x,y,z)xyw-1/2J(x,y,z)yxw.$$
 The identities  \eqref{a8} and (a) imply the required result.\\
Let $r\geq 1$ and $J(x,y,z)x_1...x_{r-1}(wt)\in J(M_3,M_3,M_3)X$. We denote $J(x,y,z)x_1...x_{r-2}$ by $J_0$. 
The identity \eqref{a1} gives
$$J_0x_{r-1}x_r(wt)=-wtJ_0x_{r-1}x_r-x_r(wt)J_0x_{r-1}-x_{r-1}x_r(wt)J_0+J_0x_r(x_{r-1}(wt))=$$
$$-wtJ_0x_{r-1}x_r+wtx_rJ_0x_{r-1}-x_{r-1}x_r(wt)J_0-J_0x_r(wtx_{r-1}). $$ 
This equality and the induction hypothesis imply  the required. $\blacksquare$ \\

{\bf Corollary 2.} Let $f,g,h$ be the polynomials from the subalgebra,
generated by elements $x$ and $y$ of the Malcev algebra and $deg f=n$. Then \\
1) $J(yx,f,z)=\frac{(-1)^n-1}{2}fJ(x,y,z),$\\
2) $J(x,f,yz)=J(x,f,z)y-((-1)^n+1)fJ(x,y,z),$\\
3) if $f$ has even degree, then $J(z,g,h)L_{f,y}=0.$\\

 {\sc Proof}. 1) Applying the operator $\Delta_y^1(f)$ to the identity $J(yx,y,z)=J(x,y,z)y$
we obtain
$$J(yx,f,z)=J(y,fx,z)+J(x,f,z)y-fJ(x,y,z).$$ 
From the identity \eqref{a13}:
$$J(yx,f,z)=-J(y,f,z)x-(-1)^nfJ(x,y,z)+J(x,f,z)y-fJ(x,y,z).$$ Using \eqref{a12} we obtain 1).\\
2) Applying the operator 
$\Delta_z^1(f)$ to identity $J(x,zy,z)=J(x,y,z)z$ we obtain $$J(x,fy,z)+J(x,zy,f)=J(x,y,z)f.$$
Using now \eqref{a13}, we obtain the required result.\\
3) We first show that $J(x,y,z)L_{f,y}=0.$ Combining \eqref{a12} and \eqref{a13}
we obtain $$J(x,f,z)y=J(fx,y,z)+\frac{(-1)^n+1}{2}fJ(x,y,z).$$
Applying this equality and  2) of this corollary, we obtain
$$J(x,y,z)fy=-J(x,f,z)yy+J(fx,y,z)y=-J(x,f,yz)y-2fJ(x,y,z)y+J(fx,y,yz).$$
That is $$J(x,y,z)fy=J(x,f,yz)y-J(fx,y,yz).$$ 
Using the identities \eqref{a12} and \eqref{a13} we transform the right part of the last equality:
$$J(x,y,z)fy=J(y,f,yz)x+2fJ(x,y,yz)+J(y,fx,yz)=$$
$$-fJ(x,y,yz)+2fJ(x,y,yz)=-J(x,y,z)yf.$$ Then, from Proposition 1 it follows that $J(z, g, h)$ can
be represented as
$$\sum_{i}\alpha_iJ(z,x,y)x_{i1}x_{i2}...x_{ik_i},$$ where
$k_i\in \mathbb{N}\cup \{0\}, \alpha_i \in F$. From  homogeneity
it follows that $x_{ij}\in \{x,y\}$. Thus, it is clear that
$$J(z,g,h)=J(\pm \sum_{i}\alpha_izx_{i1}x_{i2}...x_{ik_i},x,y).$$
This implies the required assertion. $\blacksquare$\\

We recall that the Lie center $Z(M)$ of Malcev algebra $M$ is the set of elements $c$ of algebra $M$ such that 
$J(c,a,b)=0$ for any $a$ and $b$ from the algebra $M$. It is well known (see, for example, ~\cite{sag} ) that  
the Lie center of any Malcev algebra is ideal and that $fg$ belongs to $Z(M)$ if  $J(f,g,a)=0$ for any $a$.
The central elements play important role in the theory of Malcev algebras
The  Lemma 3 ~\cite[\S 3]{fil} implies that the elements $yx^{2k-1}(yx)$ are central for any $k>0$.
The following assertion gives some central elements of more general form.

{\bf Corollary 3.} If  $f$ and $h$ are the polynomials of even degree from
the subalgebra  generated by the elements $x$ and $y$ of the Malcev algebra $M$, then  $fh$ belongs to
$Z(M).$\\

{\sc Proof.} Induction on the degree $k$ of the polynomial $h$. When
$k = 2$ the assertion obviously follows from  1) of Corollary 2.
Now suppose that the statement is true in the case when the degree
 is equal to $k = 2l$. Prove it for a polynomial $h_1$ of polynomial degree $k = 2l +2$.\\
1) Let $h_1=hxy$. We rewrite the identity \eqref{a9} as
$$J(z,x,txy)=1/2J(x,t,y)xz-1/2J(x,t,y)zx-J(x,z,t)L_{x,y}.$$
Apply the operator $\Delta_x^1(f)$ to this equation and replace  $t$
by $h$
$$J(z,f,hxy)+J(z,x,hfy)=1/2J(f,h,y)xz+1/2J(x,h,y)fz-$$
$$-1/2J(f,h,y)zx-1/2J(x,h,y)zf-J(f,z,h)L_{x,y}-J(x,z,h)L_{f,y}.$$
That is,  $J(z,f,hxy)=J(z,x,h)L_{f,x}=0.$\\
 2) Let $h_1=hxx$. Apply the operator $\Delta_x^1(f)$ to the equation $J(z,x,txx)=J(z,x,t)xx$ 
 and replace $t$ by $h$, obtain
 $$J(z,f,hxx)+J(z,x,hfx)+J(z,x,hxf)=J(z,f,h)xx+J(z,x,h)fx+J(z,x,h)xf.$$
 Hence, $J(z,f,hxx)=J(z,x,h)L_{f,x}=0.$\\
 The proofs of the cases $h_1=hyy$ and $h_1=hyx$ are similar.
$\blacksquare$\\

{\bf Lemma 3.} In any Malcev algebra $M$ the following identities hold:\\
\begin{gather}
J(aca^{2k+1}(R_cR_a)^{n-1},b,c)=J(a,b,c)L_{a,a}^{k}L_{a,c}^{n}, \label{a14} \\
J(aca^{2k+1},b,c)=J(a,b,c)L_{a,a}^{k}L_{a,c}, k\geq 0. \tag{\ref{a14}$'$}    \\        
J(aca(R_cR_a)^{n-1},b,c)=J(a,b,c)L_{a,c}^{n}, n\geq 1. \tag{\ref{a14}$''$}
\end{gather}

{\sc Proof.} We write  the identity (6) of ~\cite[\S 1]{fil}\\
$$J(ca^{2k+2},b,c)=J(a,b,c)a^{2k}(ca)-J(a,b,c)a^{2k+1}c.$$ Applying \eqref{a8} we obtain
$$J(aca^{2k+1},b,c)=-1/2J(a,ba^{2k},c)ac+1/2J(a,ba^{2k},c)ca+J(a,b,c)a^{2k+1}c.$$
That is $J(aca^{2k+1},b,c)=J(a,b,c)L_{a,a}^{k}L_{a,c}$.  Using the identities \eqref{a6} and \eqref{a9}:
$$J(a,b,c)L_{a,a}^{k}L_{a,c}^{n}=J(aca^{2k+1},b,c)L_{a,c}^{n-1}=J(aca^{2k+1}(R_cR_a)^{n-1},b,c).$$
$\blacksquare$\\

We denote $L_{zy,zy}+L_{y,y}L_{z,z}-L_{y,z}^2$ by $D(z,y).$\\

{\bf Proposition 2.} Let
$T=\{L_{x,x},L_{y,y},L_{z,z},L_{x,y},L_{x,z},L_{y,z},L_{x,zy}\}.$
For all  $S_i, T_i$ from $T$
  and for any $n$ from $\mathbb{N}\cup \{0\}$:\\
\begin{equation}\label{1prop2}
J(x,y,z)T_1T_2...T_n[S_1,S_2]=0,
\end{equation}

\begin{equation}\label{2prop2}
    J(x,y,z)T_1T_2...T_nL_{y,zy}=0,
\end{equation}

\begin{equation}\label{3prop2}
   J(x,y,z)T_1T_2...T_nD(y,z)=0. 
 \end{equation}

 {\sc Proof.} 
In  first we prove the  equalities \eqref{1prop2}, \eqref{2prop2} and \eqref{3prop2} with the assumption that $S_i$ and $T_i$ belong
 to the set $T_0=T\backslash \{L_{x,zy}\}.$ \\
In what follows we will apply the  condition $charF=0$ to collect  similar terms arising after linearization.

{\bf We prove the  equality \eqref{1prop2} with the assumption that $S_i$ and $T_i$ belong to the set $T_0=T\backslash \{L_{x,zy}\}.$}\\
1. We prove the identity

\begin{equation}
 J(x,y,z)L_{x,x}^kL_{y,y}^lL_{z,z}^m[L_{x,x},L_{y,y}]=0. \tag{a}
\end{equation}
1) $k=l=m=0$. Follows from the identity \eqref{a11}.\\
2) $k=l=0$, $m\ne 0$. 
Applying  \eqref{a6}, \eqref{a10} and again \eqref{a6} we have:

$$J(x,y,z)L_{z,z}^mL_{x,x}L_{y,y}=J(x,y,z)L_{z,z}^mxxL_{y,y}=J(x,y,xz^{2m+1})xL_{y,y}=$$
$$J(x,y,xz^{2m+1})L_{y,y}x=J(x,xy^3,xz^{2m+1})$$
Applying  \eqref{a6}, \eqref{a10} and again \eqref{a6} we have:
$$ J(x,xy^3,xz^{2m+1})=J(x,xy^3,z)L_{z,z}^{m}x=J(x,xy^3,z)xL_{z,z}^{m}=$$
$$J(x,y,z)L_{y,y}xxL_{z,z}^{m}=J(x,y,z)L_{y,y}L_{x,x}L_{z,z}^{m}$$
Applying  \eqref{a11}, \eqref{a6},  \eqref{a10} and again \eqref{a6} we have:
$$J(x,y,z)L_{y,y}L_{x,x}L_{z,z}^{m}=J(x,y,z)L_{x,x}L_{y,y}L_{z,z}^{m}=J(x,y,z)L_{x,x}yyL_{z,z}^{m}=$$ 
$$J(yx^3,y,z)yL_{z,z}^{m}=J(yx^3,y,z)L_{z,z}^{m}y=J(yx^3,y,yz^{2m+1})$$
Applying  \eqref{a6}, \eqref{a10} and again \eqref{a6} we have:
$$J(yx^3,y,yz^{2m+1})=J(x,y,yz^{2m+1})L_{x,x}y=J(x,y,yz^{2m+1})yL_{x,x}= $$
$$J(x,y,z)L_{z,z}yyL_{x,x}=J(x,y,z)L_{z,z}^mL_{y,y}L_{x,x}.$$\\
3) $k=m=0$, $l\ne 0$. It is obvious from \eqref{a11}.\\
4) $k=0$, $l\ne 0$, $m\ne 0$. From the identity \eqref{a11}:
$$J(x,y,z)L_{y,y}^lL_{z,z}^m[L_{x,x},L_{y,y}]=-J(x,y,yz^{2m+1}y^{2l-1})[L_{x,x},L_{y,y}]=0.$$
5) $k\ne 0$, $l=m=0$. Follows from the identity \eqref{a11}.\\
6) $k\ne 0$, $l=0$, $m\ne 0$. It is similar to the case 4).\\
7) $k\ne 0$, $l\ne 0$, $m=0$. Follows from \eqref{a11}.\\
8) $k\ne 0$, $l\ne 0$, $m\ne 0$. It is easy to see
$$J(x,y,z)L_{x,x}^kL_{y,y}^lL_{z,z}^m[L_{x,x},L_{y,y}]=$$
$$J(x,y,z)L_{x,x}^kL_{y,y}^lL_{z,z}^mL_{x,x}L_{y,y}-J(x,y,z)L_{x,x}^kL_{y,y}^lL_{z,z}^mL_{y,y}L_{x,x}$$
Using \eqref{a6}, \eqref{a10} and \eqref{a11} we transform the first term
$$J(x,y,z)L_{x,x}^kL_{y,y}^lL_{z,z}^mL_{x,x}L_{y,y}=-J(x,xy^{2l+1}x^{2k-1},z)L_{z,z}^mL_{x,x}L_{y,y}=$$ 
$$=J(x,xy^{2l+1}x^{2k},xz^{2m+1})L_{y,y}=J(x,xy^{2l+1},xz^{2m+1}x^{2k})L_{y,y}=$$
$$J(x,y,xz^{2m+1}x^{2k})L_{y,y}^lxL_{y,y}=J(x,y,xz^{2m+1}x^{2k})L_{y,y}^{l+1}x=$$
$$J(x,xy^{2l+3}x^{2k},xz^{2m+1})=J(x,xy^{2l+3}x^{2k},z)L_{z,z}^mx=$$
$$-J(x,xy^{2l+3}x^{2k+1},z)L_{z,z}^m$$
In the same way we can transform the second term of the commutator
$$-J(x,y,z)L_{x,x}^kL_{y,y}^lL_{z,z}^mL_{y,y}L_{x,x}=J(yx^{2k+3}y^{2l+1},y,z)L_{z,z}^m$$
Applying \eqref{a6} to the sum of the two last equalities we obtain the required. 

2. We show  that 
\begin{equation}
 J(x,y,z)L_{x,x}^kL_{y,y}^lL_{z,z}^m[L_{x,y},L_{x,x}]=0.  \tag{b}
\end{equation}
1) $k=l=m=0$. It follows from \eqref{a11}  by linearization. 
2) $m=0$.
$$J(x,y,z)L_{x,x}^kL_{y,y}^l[L_{x,y},L_{x,x}]=J(x,y,zx^{2k}y^{2l})[L_{x,y},L_{x,x}]=0.$$
3) $m\ne 0$, $k=l=0$. From the identity (a):
$$J(x,y,z)L_{z,z}^m[L_{x,x},L_{y,y}]=0.$$ Using the operator $\Delta_x^1(y)$ we obtain
$$J(x,y,z)L_{z,z}^m[L_{x,y},L_{x,x}]=0.$$
4) $m\ne 0$ and $k+l>0$. Without loss of generality we can assume that, for example, $k\ne 0$. From the identitie  (a) and \eqref{a6}:

$$J(x,y,z)L_{x,x}^kL_{y,y}^lL_{z,z}^m[L_{x,y},L_{x,x}]=-J(x,y,xz^{2m+1}x^{2k-1}y^{2l})[L_{x,y},L_{x,x}]=0.$$\\ 

3. Prove now that
\begin{equation}
 J(x,y,z)L_{x,x}^kL_{y,y}^lL_{z,z}^m[L_{x,y},L_{z,z}]=0.\tag{c}
\end{equation}
1) $k=l=m=0$. From the identity (b) follows
$J(x,y,z)[L_{x,z},L_{z,z}]=0.$ Applying the operator $\Delta_z^1(y)$ we obtain
 $$J(x,y,z)[L_{x,y},L_{z,z}]+2J(x,y,z)[L_{x,z},L_{y,z}]=0.$$
Appplying the identity (\ref{a14}$''$) to each term of the commutator we transform
$$J(x,y,z)[L_{x,z},L_{y,z}]=J(x,y,z)L_{x,z}L_{y,z}-J(x,y,z)L_{y,z}L_{x,z}=J(xzx,y,z)L_{y,z}-$$
$$-J(x,yzy,z)L_{x,z}=J(xzx,yzy,z)-J(xzx,yzy,z)=0.$$
Hence,  $J(x,y,z)[L_{x,y},L_{z,z}]=0$.\\
 2) $k=l=0$, $m\ne 0$. From the identity (a) follows $J(x,y,z)L_{z,z}^m[L_{x,x},L_{z,z}]=0.$
 Applying the operator $\Delta_x^1(y)$  we obtain $J(x,y,z)L_{z,z}^m[L_{x,y},L_{z,z}]=0.$\\
 3) $k=0$, $l\ne 0$, $m=0$. From the identity (a) it follows $J(x,y,z)L_{y,y}^l[L_{x,x},L_{z,z}]=0.$
 Applying the operator $\Delta_x^1(y)$ we have $J(x,y,z)L_{y,y}^l[L_{x,y},L_{z,z}]=0.$\\
 4) $k=0$, $l\ne 0$, $m\ne 0$. From the identity (a) it follows $J(x,y,z)L_{y,y}^lL_{z,z}^m[L_{x,x},L_{z,z}]=0.$
 Applying the operator  $\Delta_x^1(y)$  we obtain $J(x,y,z)L_{y,y}^lL_{z,z}^m[L_{x,y},L_{z,z}]=0.$\\
 5) $k\ne 0$, $l=m=0$. From the identity (a) it follows $J(x,y,z)L_{x,x}^k[L_{y,y},L_{z,z}]=0.$
Applying the operator $\Delta_y^1(x)$  we obtain
$J(x,y,z)L_{x,x}^k[L_{x,y},L_{z,z}]=0.$\\
6) $k\ne 0$, $l=0$, $m\ne 0$. From the identity (a) follows
$J(x,y,z)L_{x,x}^kL_{z,z}^m[L_{y,y},L_{z,z}]=0.$ Applying the operator
$\Delta_y^1(x)$ we have
$J(x,y,z)L_{x,x}^kL_{z,z}^m[L_{x,y},L_{z,z}]=0.$\\
7) $k\ne 0$, $l\ne 0$, $m=0$. 
$$J(x,y,z)L_{x,x}^kL_{y,y}^l[L_{x,y},L_{z,z}]=$$
$$J(x,y,z)L_{x,x}^kL_{y,y}^lL_{x,y}L_{z,z}-J(x,y,z)L_{x,x}^kL_{y,y}^lL_{z,z}L_{x,y}$$
Applying the identities (b) and (\ref{a14}$'$) to the first term of commutator we obtain

$$J(x,y,z)L_{x,x}^kL_{y,y}^lL_{x,y}L_{z,z}=J(x,y,z)L_{x,x}^kL_{x,y}L_{y,y}^lL_{z,z}=J(xyx^{2k+1},y,z)L_{y,y}^lL_{z,z}$$

From \eqref{a11} and \eqref{a6}
$$J(xyx^{2k+1},y,z)L_{y,y}^lL_{z,z}=J(xyx^{2k+1},y,z)L_{z,z}L_{y,y}^l=J(xyx^{2k+1},y,z)L_{z,z}yyL_{y,y}^{l-1}=$$
$$J(xyx^{2k+1},y,yz^3)yL_{y,y}^{l-1}=J(xyx^{2k+1},y,yz^3)y^{2l-1}=-J(xyx^{2k+1},y,yz^3y^{2l-1})$$

From (\ref{a14}$'$), \eqref{a6} and (a)

$$-J(xyx^{2k+1},y,yz^3y^{2l-1})=-J(x,y,yz^3y^{2l-1})L_{x,x}^kL_{x,y}=J(x,y,yz^3)yL_{y,y}^{l-1}L_{x,x}^kL_{x,y}=$$
$$J(x,y,z)L_{y,y}^lL_{z,z}L_{x,x}^kL_{x,y}$$
In the same way we can transform the second term of the commutator
$$-J(x,y,z)L_{x,x}^kL_{y,y}^lL_{z,z}L_{x,y}=-J(x,y,z)L_{x,x}^kL_{y,y}^lL_{z,z}L_{x,y}$$
Applying (a) to the sum of the two last equalities we obtain the required. 

8) $k\ne 0$, $l\ne 0$, $m\ne 0$. Induction on $m$. For $m=0$  this statement was proved in 7). 
We suppose that 
$$ J(x,y,z)L_{x,x}^kL_{y,y}^lL_{z,z}^i[L_{x,y},L_{z,z}]=0$$ for all $i<m$. 
$$J(x,y,z)L_{x,x}^kL_{y,y}^lL_{z,z}^m[L_{x,y},L_{z,z}]=$$
$$J(x,y,z)L_{x,x}^kL_{y,y}^lL_{z,z}^mL_{x,y}L_{z,z}-J(x,y,z)L_{x,x}^kL_{y,y}^lL_{z,z}^{m+1}L_{x,y}$$ 
We transform the first term of the commutator. From the induction hypothesis, (a), (\ref{a14}$'$) and \eqref{a6}  we obtain
$$J(x,y,z)L_{x,x}^kL_{y,y}^lL_{z,z}^mL_{x,y}L_{z,z}=J(x,y,z)L_{x,x}^kL_{y,y}^lL_{x,y}L_{z,z}^{m+1}=$$
$$J(x,y,z)L_{x,x}^kL_{x,y}L_{y,y}^lL_{z,z}^{m+1}=J(yx^{2k+1},y,z)L_{y,y}^lL_{z,z}^{m+1}=$$
$$J(yx^{2k+1},y,z)L_{z,z}^{m+1}L_{y,y}^l=-J(yx^{2k+1},y,yz^{2m+3}y^{2l-1})=$$
$$-J(x,y,yz^{2m+3}y^{2l-1})L_{x,x}^kL_{x,y}=J(x,y,z)L_{y,y}^lL_{z,z}^{m+1}L_{x,x}^kL_{x,y}=$$
$$J(x,y,z)L_{x,x}^kL_{y,y}^lL_{z,z}^{m+1}L_{x,y}$$

4. We now prove  the identity
\begin{equation}
 J(x,y,z)L_{x,x}^kL_{y,y}^lL_{z,z}^m[L_{x,y},L_{x,z}]=0 \tag{d}
\end{equation}
1) $m=0.$ From  (c)
using the operator $\Delta_z^1(x).$\\
2) $l=0$. From the identity (c):
$$J(x,y,z)L_{x,x}^kL_{z,z}^m[L_{y,y},L_{x,z}]=0.$$ Applying the operator $\Delta_y^1(x),$  we obtain the required result.\\
3) $k=0$. The proof is similar to the cases  1) and 2).\\
4) $k\ne 0$, $l\ne 0$, $m\ne 0$. Using the identities (a), (c), (b) and (\ref{a14}$'$):
 
$$J(x,y,z)L_{x,x}^kL_{y,y}^lL_{z,z}^mL_{x,y}L_{x,z}=J(x,y,z)L_{y,y}^lL_{x,y}L_{x,x}^kL_{z,z}^mL_{x,z}=$$
$$J(x,yxy^{2l+1},z)L_{x,x}^kL_{z,z}^mL_{x,z}=J(x,yxy^{2l+1},zxz^{2m+1}x^{2k})=$$
$$J(x,y,zxz^{2m+1}x^{2k})L_{y,y}^lL_{x,y}=J(x,y,z)L_{x,x}^kL_{z,z}^mL_{x,z}L_{y,y}^lL_{x,y}=$$
$$J(x,y,z)L_{x,x}^kL_{y,y}^lL_{z,z}^mL_{x,z}L_{x,y}.$$ It finishes the proof of the equality (d).\\

From the proven identities (a),(b),(c) and (d) it follows
$$J(x,y,z)L_{x,x}^kL_{y,y}^lL_{z,z}^m[S_1,S_2]=0, \textrm{where } S_i \in T_0.$$ 
Hence, with the operator $\Delta_x^1(y)$ it is easy to prove the identity
$$J(x,y,z)L_{x,y}^nL_{x,x}^kL_{y,y}^lL_{z,z}^m[S_1,S_2]=0,$$ using induction on $n$.
The induction on $p$ and operator $\Delta_z^1(y)$  give the identity
$$J(x,y,z)L_{x,y}^nL_{y,z}^pL_{x,x}^kL_{y,y}^lL_{z,z}^m[S_1,S_2]=0.$$
Finally,  induction on $q$ and the operator $\Delta_x^1(z)$ give the identity
$$J(x,y,z)L_{x,y}^nL_{y,z}^pL_{x,z}^qL_{x,x}^kL_{y,y}^lL_{z,z}^m[S_1,S_2]=0.$$
This  is the equality \eqref{1prop2} assuming that the $ S_i $ and $ T_i $
belong to the 
 set $T_0=T\backslash \{L_{x,zy}\}.$ \\

{\bf We prove the equality  \eqref{2prop2} with the assumption that $T_i$ belong to the set $T_0=T\backslash \{L_{x,zy}\}.$}\\
That is, 
\begin{equation} 
J(x,y,z)L_{x,y}^nL_{x,z}^qL_{y,z}^pL_{x,x}^kL_{y,y}^lL_{z,z}^mL_{y,zy}=0. \tag{e} 
\end{equation} 
Induction on q. In first we consider the case $q=0.$ That is 
\begin{equation} 
  J(x,y,z)L_{x,y}^nL_{y,z}^pL_{x,x}^kL_{y,y}^lL_{z,z}^mL_{y,zy}=0 \tag{e$'$}
\end{equation}
 1) $n=k=0.$  Obviously follows from the Proposition 1 and from equality 3) of the Corollary 2.\\
 2) $n=0$, $k\ne 0$.\\   
     2.a) $p=l=m=0.$ Applying the identity \eqref{a6} and identities already proven in this proposition we obtain
      $$J(x,y,z)L_{x,x}^kyL_{y,z}=J(yx^{2k+1},y,z)L_{y,z}=J(yx^{2k+1},y,zyz)=J(x,y,zyz)L_{x,x}^ky=$$
      $$J(x,y,z)L_{y,z}L_{x,x}^ky=J(x,y,z)L_{x,x}^kL_{y,z}y.$$ That is,
      $J(x,y,z)L_{x,x}^kL_{y,z}y-J(x,y,z)L_{x,x}^kyL_{y,z}=0.$ From the identity \eqref{a7}:
      $$0=J(x,y,z)L_{x,x}^kL_{y,z}y-J(x,y,z)L_{x,x}^kyL_{y,z}=J(x,y,z)L_{x,x}^kL_{y,yz}.$$
     2.b) $p\ne 0,$ $l=m=0$
  Using the identity \eqref{a6} we obtain
   $$J(x,y,z)L_{y,z}^pL_{x,x}^kL_{y,zy}=1/2J(x,y,z)L_{x,x}^kyzL_{y,z}^{p-1}L_{y,zy}+1/2J(x,y,z)L_{x,x}^kzyL_{y,z}^{p-1}L_{y,zy}=$$
   $$-1/2J(yx^{2k+1}L_{y,z}^{p-1}z,y,z)L_{y,zy}-1/2J(zx^{2k+1}L_{y,z}^{p-1}y,y,z)L_{y,zy}=0.$$
   2.c) If $p=0$ and $l\ne 0$ or $p=0$ and $m\ne 0$ it is sufficient to apply the identity \eqref{a6}
 3) $n\ne 0$, $k\ne 0$. Using identity \eqref{a14} 
 $$J(x,y,z)L_{x,y}^nL_{y,z}^pL_{x,x}^kL_{y,y}^lL_{z,z}^mL_{y,zy}=J(x,y,z)L_{x,x}^kL_{x,y}^nL_{y,z}^pL_{y,y}^lL_{z,z}^mL_{y,zy}=$$
 $$J(xyx^{2k+1}(R_xR_y)^{n-1},y,z)L_{y,z}^pL_{y,y}^lL_{z,z}^mL_{y,zy}=$$
 $$J(xyx^{2k+1}(R_xR_y)^{n-1}L_{y,z}^pL_{y,y}^lL_{z,z}^m,y,z)L_{y,zy}=0$$  
 4) $n\ne 0$, $k=0$. It is sufficient to apply (\ref{a14}$''$).\\
Applying the operator  $\Delta_x^1(z)$ to the identity (e$'$) and using 
induction hypothesis  it is easy to prove the required identity (e).
\\

{\bf We prove the equality \eqref{3prop2} with the assumption that $T_i$ belong to the set $T_0=T\backslash \{L_{x,zy}\}.$}\\

 1) $n=p=k=l=m=0$.  
We transform the first term of this sum. Substituing in the identity \eqref{a8} $t$ by $b$ we obtain
$J(a,b,c)=-1/2J(a,b,c)ab+1/2J(a,b,c)ba.$ Applying this identity to the first term of the sum
$$J(x,y,z)D(y,z)=J(x,y,z)L_{yz,yz}+J(x,y,z)L_{y,y}L_{z,z}-J(x,y,z)L_{y,z}^2$$  we have:
$$J(x,y,z)(yz)(yz)=1/2J(-xyz+xzy,y,z)(yz)=$$
$$ -1/4J(-xyz+xzy,y,z)yz+1/4J(-xyz+xzy,y,z)zy=$$
$$1/4(J(x,y,z)yzyz-J(x,y,z)yzzy-J(x,y,z)zyyz+J(x,y,z)zyzy)=$$
$$-1/4J(x,y,z)yyzz+1/2J(x,y,z)yL_{y,z}z-1/4J(x,y,z)yL_{z,z}y-1/4J(x,y,z)zL_{y,y}z-$$
$$-1/4J(x,y,z)zzyy+1/2J(x,y,z)zL_{y,z}y=-J(x,y,z)L_{y,y}L_{z,z}+J(x,y,z)L_{y,z}^2$$
That is,  $J(x,y,z)D(y,z)=0.$\\
 2) $n=p=l=m=0$, $k\ne 0.$
Apply to the identity
$$J(x,y,z)L_{x,x}^kL_{y,y}=J(x,y,z)L_{y,y}L_{x,x}^k$$
the operator $\Delta_y^2(zy)$:
$$J(x,y,z)L_{x,x}^kL_{zy,zy}+2J(x,zy,z)L_{x,x}^kL_{y,zy}=J(x,y,z)L_{zy,zy}L_{x,x}^k+2J(x,zy,z)L_{y,zy}L_{x,x}^k,$$
We show that the second term in each part of this equality is zero.  
Really, from the identity \eqref{a6} it follows
$J(x,zy,z)L_{x,x}^kL_{y,zy}=J(x,y,z)zL_{x,x}^kL_{y,zy}=J(zx^{2k+1},y,z)L_{y,zy}=0.$ And $J(x,zy,z)L_{y,zy}L_{x,x}^k=0.$
Therefore, 
$J(x,y,z)L_{x,x}^kL_{zy,zy}=J(x,y,z)L_{zy,zy}L_{x,x}^k.$
From the  identities proven above and item 1) it follows
$$J(x,y,z)L_{x,x}^kD(y,z)=J(x,y,z)D(y,z)L_{x,x}^k=0,$$ that
required. 
The remaining cases are proved similarly to the equality \eqref{2prop2}.\\

{\bf We prove the equality \eqref{1prop2} with the assumption that $T_i$ belong to the set $T_0$ and $S_j$ belongs to the set $T$.}\\
That is, \begin{equation}
          J(x,y,z)T_1T_2...T_n[S_1,S_2]=0,\quad T_i \in T_0, \quad S_j \in T   \tag{f}
         \end{equation}

Let none of the operators $ T_i $ be not equal $L_{x,x}$. From the  equalities proven above it follows  
$$J(x,y,z)T_1T_2...T_n[L_{x,x},L_{y,y}]=0.$$
Apply the operator $\Delta_x^1(zy).$ Obtain
$L_{x,y}\Delta_x^1(zy)=L_{zy,y},$ $L_{y,y}\Delta_x^1(zy)=0,$ 
$L_{z,z}\Delta_x^1(zy)=0$,\quad $L_{x,z}\Delta_x^1(zy)=L_{zy,z},$\quad
$L_{y,z}\Delta_x^1(zy)=0$ and, therefore,  
$$J(x,y,z)T_1T_2...T_n[L_{x,zy},L_{y,y}]=0, \textrm{for } T_i\ne L_{x,x}.$$
 The proof of  identity (f)  is obtained by induction on the
degree of the operator   $L_{x,x}$ using the operator
$\Delta_{z}^1(x)$.  
Really, $L_{x,z}\Delta_{z}^1(x)=L_{x,x}$ and
 the remaining terms arising from the action of this 
 operator are equal to zero by the previously proven identities. Applying to the resulting identity   corresponding operator
$\Delta_{x_i}^1(x_j), \textrm{where } x_i, x_j\in X,$ we obtain the equality (f).
\\

{\bf We prove all statements of this Proposition without  restrictions on
$T_i$ and $S_j$.}
We shall prove the conjunction of all three statements using
induction on $l$, where $l$ is number of operators $L_{x,zy}$ in
the sequence $T_1, T_2,...,T_n$. For $l=0$ all three identities are
proved. Assume now that all three identities hold for
all $l\leq k$. We show that these identities
hold for $l=k+1$.

First, we note that the induction hypothesis for \eqref{1prop2} implies that if we have  
any  sequence of operators from $T$ and it contains not more 
than $k+1$ copies of $L_{x,zy}$, then we can permute these operators acting on  $J(x,y,z)$.
 Suppose now that among the  $T_1, T_2,...,T_n$ there are $k+1$ copies of
$L_{x,zy}$. Let $T_i=L_{x,zy}$. Replace it by the operator
$L_{x,x}$.  From the induction hypothesis:
$$J(x,y,z)T_1T_2...L_{x,x}...T_n[S_1,S_2]=0.$$ Applying to this identity the operator $\Delta_{x}^1(zy)$ 
and taking into account equalities: $L_{x,x}\Delta_{x}^1(zy)=2L_{x,zy}$,
$L_{x,zy}\Delta_{x}^1(zy)=L_{zy,zy}$,
$L_{x,z}\Delta_{x}^1(zy)=2L_{z,zy}$,
$L_{x,y}\Delta_{x}^1(zy)=2L_{y,zy}$,
$L_{y,y}\Delta_{x}^1(zy)=L_{z,z}\Delta_{x}^1(zy)=L_{y,z}\Delta_{x}^1(zy)=0,$
  and equalities \eqref{2prop2} and \eqref{3prop2} for $l\leq k$ we obtain
$$J(x,y,z)T_1T_2...L_{x,zy}...T_n[S_1,S_2]=0.$$
Similar arguments prove  equalities \eqref{2prop2} and \eqref{3prop2} for a sequence
$T_1, T_2,...,T_n$ containing $k+1$ copies of
$L_{x,zy}$.
  $\blacksquare$\\

{\bf Corollary 4.} Let $w=J(x,y,z)T_1T_2...T_n$, for some
 $$T_i\in T=\{L_{x,x},L_{y,y},L_{z,z},L_{x,y},L_{x,z},L_{y,z},L_{x,zy}\},$$
  $k,l,m,n\in \mathbb{N}\cup \{0\}$. The following equalities hold:\\
\begin{equation}\label{b18}
wL_{x,zy}=wL_{y,xz}=wL_{z,yx}
\end{equation}
\begin{equation}\label{b19}
w[L_{x_1,x_2},R_{x_3}]=0, x_i\in X
\end{equation}
$$J(x,y,z)L_{x,x}^kL_{y,y}^lL_{z,z}^m(xz)L_{x,xy}=$$
\begin{equation}\label{b20}
=J(x,y,z)L_{x,x}^kL_{y,y}^lL_{z,z}^mxzL_{x,xy}=J(x,y,z)L_{x,x}^kL_{y,y}^lL_{z,z}^mzxL_{x,xy}=0.
\end{equation}
\begin{equation}\label{b21}
J(x,y,z)L_{x,x}^kL_{y,y}^lL_{z,z}^mzx[R_x,L_{x,y}]=0.
\end{equation}

{\sc Proof.} The equality \eqref{b18}. From the \eqref{2prop2}:
 $wL_{x,xy}=0$. Apply the operator $\Delta_{x}^1(z)$ to this equality. 
 From the \eqref{1prop2} we obtain
$wL_{z,xy}+wL_{x,zy}=0$. That is, $wL_{z,yx}=wL_{x,zy}$.
Similarly, the operator $\Delta_{y}^1(z)$, applied to the identity
$wL_{y,yx}=0$ gives $wL_{z,yx}=wL_{y,xz}$.\\

The equality \eqref{b19}. The equalities $w[L_{x,y},R_x]=w[L_{x,y},R_y]=0,$
$w[L_{x,z},R_x]=w[L_{x,z},R_z]=0,$
$w[L_{y,z},R_y]=w[L_{y,z},R_z]=0$ follow from the identity \eqref{a7} and equality \eqref{2prop2}.
The equalities $w[L_{x,y},R_z]=w[L_{y,z},R_x]=w[L_{x,z},R_y]=0$ follow from the identity \eqref{a7} and from the equality \eqref{b18}.\\

The equality \eqref{b20}. From the  equalities \eqref{b19} and \eqref{1prop2}:
$$J(x,y,z)L_{y,y}^lL_{z,z}^mzL_{x,x}^kL_{x,y}=J(x,y,z)L_{x,y}L_{y,y}^lL_{z,z}^mzL_{x,x}^k.$$
Apply to this identity the operator $\Delta_{y}^2(xy)$. 
From the \eqref{1prop2} and \eqref{2prop2} we have:
$$J(x,xy,z)L_{y,y}^{l-i}L_{y,xy}L_{y,y}^{i-1}L_{z,z}^mzL_{x,x}^kL_{x,y}=J(x,y,xz)L_{x,y}L_{y,y}^{l-i}L_{y,xy}L_{y,y}^{i-1}L_{z,z}^mzL_{x,x}^k=0,$$
$$J(x,y,z)L_{y,y}^{l-i}L_{y,xy}L_{y,y}^{i-1}L_{z,z}^mzL_{x,x}^kL_{x,xy}=J(x,y,z)L_{x,xy}L_{y,y}^{l-i}L_{y,xy}L_{y,y}^{i-1}L_{z,z}^mzL_{x,x}^k=0,$$
$$J(x,y,z)L_{y,y}^{l-i}L_{xy,xy}L_{y,y}^{i-1}L_{z,z}^mzL_{x,x}^kL_{x,y}=J(x,y,z)L_{x,y}L_{y,y}^{l-i}L_{xy,xy}L_{y,y}^{i-1}L_{z,z}^mzL_{x,x}^k,$$
$$J(x,xy,z)L_{x,xy}L_{y,y}^{l}L_{z,z}^mzL_{x,x}^k=0.$$
Therefore, the application of the operator gives 
$$
J(x,xy,z)L_{y,y}^lL_{z,z}^mzL_{x,x}^kL_{x,xy}=0.
$$
Furthermore,
$$0=J(x,xy,z)L_{y,y}^lL_{z,z}^mzL_{x,x}^kL_{x,xy}=J(x,y,z)L_{y,y}^lL_{z,z}^mxzL_{x,x}^kL_{x,xy}=$$
$$-J(x,y,z)L_{y,y}^lL_{z,z}^mzxL_{x,x}^kL_{x,xy}+2J(x,y,z)L_{y,y}^lL_{z,z}^mL_{x,z}L_{x,x}^kL_{x,xy}=$$
$$-J(x,y,z)L_{y,y}^lL_{z,z}^mzxL_{x,x}^kL_{x,xy}=-J(x,y,z)L_{y,y}^lL_{z,z}^mzxx^{2k}L_{x,xy}=$$
$$-J(x,y,z)L_{y,y}^lL_{z,z}^mzx^{2k}xL_{x,xy}=-J(x,y,z)L_{y,y}^lL_{z,z}^mL_{x,x}^kzxL_{x,xy}=$$
$$J(x,y,z)L_{y,y}^lL_{z,z}^mL_{x,x}^kxzL_{x,xy}-2J(x,y,z)L_{y,y}^lL_{z,z}^mL_{x,x}^kL_{x,z}L_{x,xy}=$$
$$J(x,y,z)L_{y,y}^lL_{z,z}^mL_{x,x}^kxzL_{x,xy}.$$ 
That is,
\begin{equation}
J(x,y,z)L_{y,y}^lL_{z,z}^mL_{x,x}^kxzL_{x,xy}=-J(x,y,z)L_{y,y}^lL_{z,z}^mL_{x,x}^kzxL_{x,xy}=0. \tag {a} 
\end{equation}

We prove now the equality $$J(x,y,z)L_{x,x}^kL_{y,y}^lL_{z,z}^m(xz)L_{x,xy}=0.$$
1) $m>0$.  Then 
$$J(x,y,z)L_{y,y}^lL_{z,z}^mL_{x,x}^k(xz)L_{x,xy}=-J(x,zy^{2l+1}z^{2m-1}x^{2k},z)(xz)L_{x,xy}=$$
$$-1/2J(x,zy^{2l+1}z^{2m-1}x^{2k},z)zxL_{x,xy}+1/2J(x,zy^{2l+1}z^{2m-1}x^{2k},z)xzL_{x,xy}=$$
$$1/2J(x,y,z)L_{y,y}^lL_{z,z}^mL_{x,x}^kzxL_{x,xy}-1/2J(x,y,z)L_{y,y}^lL_{z,z}^mL_{x,x}^kxzL_{x,xy}=0.$$
2) $k>0$ It is similar to the case 1).\\
3) $k=l=m=0$. Easy to see from \eqref{a8}
$$J(x,y,z)(xz)L_{x,xy}=1/2J(x,y,z)zxL_{x,xy}-1/2J(x,y,z)xzL_{x,xy}=$$
$$-J(x,y,z)xzL_{x,xy}=-J(x,xy,z)zL_{x,xy}=-J(x,xy,z)L_{x,xy}z=0$$
4) $k=m=0$, $l\ne 0$. From \eqref{a8} and \eqref{2prop2}:
$$J(x,y,z)(zx)L_{y,y}^lL_{x,xy}=1/2J(x,y,z)xzL_{y,y}^lL_{x,xy}-1/2J(x,y,z)zxL_{y,y}^lL_{x,xy}=$$
\begin{equation}
J(x,y,z)xzL_{y,y}^lL_{x,xy}.\tag{b}
\end{equation}
 Applying the operator $\Delta_{z}^1(zx)$ to the identity
$$J(x,y,z)L_{y,y}^lz=J(x,y,z)zL_{y,y}^l$$ we obtain
$$J(x,y,zx)L_{y,y}^lz+J(x,y,z)L_{y,y}^l(zx)=J(x,y,z)(zx)L_{y,y}^l+J(x,y,zx)zL_{y,y}^l$$ and
$$-J(x,y,z)L_{y,y}^lxzL_{x,xy}+J(x,y,z)L_{y,y}^l(zx)L_{x,xy}=$$
$$J(x,y,z)(zx)L_{y,y}^lL_{x,xy}-J(x,y,z)xzL_{y,y}^lL_{x,xy}.$$

From the identities (a) and (b) obtain:
$$J(x,y,z)L_{y,y}^l(zx)L_{x,xy}=0.$$ 
The identity \eqref{b21} follows from \eqref{b20} and \eqref{a7}.
$\blacksquare$\\

{\bf Lemma 4.} In any Malcev algebra $M$ the following identities hold:\\
\begin{equation}\label{b24}
 t(yx)z+t(zy)x+t(xz)y=tL_{x,zy}+tL_{y,xz}+tL_{z,yx}-1/2J(x,y,z)t,
\end{equation}
$$G(t,x,y,z)+J(x,y,z)t=2/3(tyzx-txzy)+2/3(tzxy-tyxz)+$$
\begin{equation}\label{b25}
+2/3(txyz-tzyx)+2/3tL_{x,zy}+2/3tL_{y,xz}+2/3tL_{z,yx}.
\end{equation}
\begin{equation}\label{b22}
G(tx,x,y,z)=G(t,x,y,z)x-2J(x,y,z)(tx),
\end{equation}
\begin{equation}\label{b23}
G(tx^2,x,y,z)=G(t,x,y,z)x^2+2J(x,t,z)L_{x,xy}+2J(x,y,t)L_{x,xz},
\end{equation}

 {\sc Proof.}
The equality \eqref{b24}.
 From the identity \eqref{a1} we have:
$$t(xy)z+t(yz)x+t(zx)y=zxyt+tzxy-tyzx-xt(yz)+xyzt+txyz-tzxy-yt(zx)+$$
$$yzxt+tyzx-txyz-zt(xy)=ty(zx)+tx(yz)+ty(zx)+J(x,y,z)t.$$
 That is,$$t(xy)z+t(yz)x+t(zx)y=ty(zx)+tx(yz)+tz(xy)+J(x,y,z)t.$$
Hence   \eqref{b24} follows.\\

The equality \eqref{b25}.
 Further, from this identity and the identities \eqref{a1} and \eqref{a3} we obtain:
$$3/2G(t,x,y,z)+3/2J(x,y,z)t=J(t,y,z)x+J(x,t,z)y+J(x,y,t)z-J(x,y,z)t=$$
$$tyzx+ztyx+yztx+xtzy+zxty+tzxy+xytz+txyz+ytxz=$$
$$tyzx-txzy+tzxy-tyxz+txyz-tzyx-t(xy)z+t(yz)x+t(zx)y=$$
$$tyzx-txzy+tzxy-tyxz+txyz-tzyx+tL_{x,zy}+tL_{y,xz}+tL_{z,yx}.$$

 The equality \eqref{b22}. From the identities
\eqref{a3}, \eqref{a4} and the definition of the function $G$:
$$3/2G(tx,x,y,z)=J(tx,y,z)x+J(x,tx,z)y+J(x,y,tx)z-J(x,y,z)(tx)=$$
$$J(tx,y,z)x+J(yz,x,tx)-G(y,z,x,tx)-J(x,y,z)(tx)=$$
$$J(tx,y,z)x-J(yz,x,t)x+G(tx,x,y,z)-J(x,y,z)(tx)=$$
$$G(tx,x,y,z)-J(x,y,z)(tx)+1/2G(t,x,y,z)x,$$ 
whence follows the desired identity.\\
The identity \eqref{b23}. From the identities \eqref{a8}, \eqref{a3}, \eqref{a4} and the definition of the function
$G$:
$$(tx)J(x,y,z)=-1/2J(x,t,z)xy-1/2J(t,x,y)xz-J(y,z,tx)x-3/2J(x,t,zy)x=$$
$$1/2J(x,t,z)yx+1/2J(t,x,y)zx-J(y,z,tx)x-3/2J(x,t,zy)x-J(x,t,z)L_{x,y}-J(t,x,y)L_{x,z}=$$
$$1/2(J(x,t,z)y+J(t,x,y)z+J(t,x,zy))x-1/2G(t,x,y,z)x-J(x,t,z)L_{x,y}-J(t,x,y)L_{x,z}=$$
$$1/2(G(y,z,t,x)-2J(t,x,yz))x-1/2G(t,x,y,z)x-J(x,t,z)L_{x,y}-J(t,x,y)L_{x,z}=$$
$$-J(x,t,z)L_{x,y}-J(t,x,y)L_{x,z}-J(t,x,yz)x.$$ 
Further, from the identity \eqref{b22} we have:
$$G(tx,x,y,z)=G(t,x,y,z)x-2J(x,y,z)(tx)=$$
$$G(t,x,y,z)x-2J(x,t,z)L_{x,y}-2J(x,y,t)L_{x,z}-2J(t,x,yz)x.$$
Thus, from this identity and the identity \eqref{a7} we obtain:
$$G(tx^2,x,y,z)=G(tx,x,y,z)x-2J(x,tx,z)L_{x,y}-2J(tx,x,y)L_{x,z}-2J(tx,x,yz)x=$$
$$G(t,x,y,z)x^2-2J(x,t,z)L_{x,y}x-2J(x,y,t)L_{x,z}x-2J(t,x,yz)x^2+2J(x,t,z)xL_{x,y}+$$
$$+2J(t,x,y)xL_{x,z}+2J(t,x,yz)x^2=G(t,x,y,z)x^2+2J(x,t,z)L_{x,xy}+2J(x,y,t)L_{x,xz}.$$
It implies \eqref{b23}. $\blacksquare$ \\

{\bf Lemma 5.} Let
$T=\{L_{x,x},L_{y,y},L_{z,z},L_{x,y},L_{x,z},L_{y,z},L_{x,zy}\},$
and $T_i \in T$, $n\in \mathbb{N}\cup \{0\}.$
 In the algebra $M_3$ the following relations hold:\\
\begin{equation}\label{b26}
J(x,y,z)T_1T_2...T_nJ(x,y,z)=0,
\end{equation}
\begin{equation}\label{b27}
J(x,y,z)T_1T_2...T_n(J(x,y,z)L_{x,x})=0,
\end{equation}
\begin{equation}\label{b28}
J(x,y,z)T_1T_2...T_nx(J(x,y,z)x)=0.
\end{equation}

{\sc Proof.} We prove \eqref{b26}.
We first prove that
$$J(x,y,z)L_{x,x}^kL_{y,y}^lL_{z,z}^mJ(x,y,z)=0.$$ 
Further, we denote $J(x,y,z)L_{x,x}^kL_{y,y}^lL_{z,z}^m$ by  $v$.
From the identity \eqref{b23} :
$$G(vL_{x,x},x,y,z)-G(v,x,y,z)L_{x,x}=2J(x,v,z)L_{x,xy}+2J(x,y,v)L_{x,xz}.$$
If $k+m\ne 0$ and $k+l\ne 0$, then from identity \eqref{b20} and \eqref{a5} we obtain:
$$G(vL_{x,x},x,y,z)-G(v,x,y,z)L_{x,x}=$$
$$2J(x,J(x,y,z)L_{x,x}^kL_{y,y}^lL_{z,z}^m,z)L_{x,xy}+2J(x,y,J(x,y,z)L_{x,x}^kL_{y,y}^lL_{z,z}^m)L_{x,xz}=$$
$$-2J(x,J(x,xy^{2l+1}x^{2k-1}z^{2m},z),z)L_{x,xy}-2J(x,y,J(x,y,xz^{2m+1}x^{2k-1}y^{2l}))L_{x,xz}=$$
$$6J(x,y,z)L_{x,x}^kL_{y,y}^lL_{z,z}^m(xz)L_{x,xy}+6J(x,y,z)L_{x,x}^kL_{y,y}^lL_{z,z}^m(xy)L_{x,xz}=0.$$
Suppose for example $k+l=0$, that is $v=J(x,y,z)L_{z,z}^m.$ Easy to see
$$J(x,y,J(x,y,z)L_{y,y}^m)L_{x,xz}=J(x,y,J(x,y,zy^{2m}))L_{x,xz}=3J(x,y,z)L_{y,y}^m(xy)L_{x,xz}=0.$$
That is, $J(x,y,J(x,y,z)L_{y,y}^m)L_{x,xz}=0.$ Applying the operator $\Delta_{y}^{2m}(z)$ to this
identity we obtain
$J(x,y,J(x,y,z)L_{z,z}^m)L_{x,xz}=0$. Further,
$$G(vL_{x,x},x,y,z)-G(v,x,y,z)L_{x,x}=2J(x,v,z)L_{x,xy}+2J(x,y,v)L_{x,xz}=0.$$
Thus, we have $G(vL_{x,x},x,y,z)-G(v,x,y,z)L_{x,x}=0,$
for all $t$ of the form $v=J(x,y,z)L_{x,x}^kL_{y,y}^lL_{z,z}^m$.
Using this equality and induction on $k+l+m$ we obtain:
\begin{equation}
J(x,y,z)L_{x,x}^kL_{y,y}^lL_{z,z}^mG=J(x,y,z)GL_{x,x}^kL_{y,y}^lL_{z,z}^m.\tag{a}
\end{equation}
We have:
$$vL_{z,z}xyz=vzzxyz=-vzxzyz+2vzL_{x,z}yz=vzxyzz-2vzxL_{y,z}z+2vzL_{x,z}yz=$$
$$-vxzyzz+2vL_{x,z}yzz-2vzxL_{y,z}z+2vzL_{x,z}yz=$$
$$vxyzzz-2vxL_{y,z}zz+2vL_{x,z}yzz-2vzxL_{y,z}z+2vzL_{x,z}yz=$$
$$vxyzL_{z,z}-2vxL_{z,z}L_{y,z}+2vL_{x,z}yL_{z,z}-2vzxL_{y,z}z+2vzL_{x,z}yz=$$
$$vxyzL_{z,z}-2vxL_{z,z}L_{y,z}+2vL_{x,z}yL_{z,z}+2vxzL_{y,z}z-2vL_{x,z}zzy-$$
$$-4vL_{x,z}L_{y,z}z+4vzL_{x,z}L_{y,z}=$$
$$vxyzL_{z,z}-2vxL_{z,z}L_{y,z}+2vxzL_{y,z}z.$$
From \eqref{b21} it follows that the sum of the last two terms is 0. That is, 
$$vL_{z,z}xyz=vxyzL_{z,z}.$$ 
Similarly,  the equalities
$vL_{y,y}xyz=vxyzL_{y,y}$ and $vL_{x,x}xyz=vxyzL_{x,x}$ can be verified.
Therefore,
\begin{equation}
J(x,y,z)L_{x,x}^kL_{y,y}^lL_{z,z}^mxyz=J(x,y,z)xyzL_{x,x}^kL_{y,y}^lL_{z,z}^m.\tag{b}
\end{equation}
Recall that we agreed to denote the action of the operator $t\mapsto G(t,x,y,z)$ by $tG$ and the superposition of $n$ these 
operators by $tG^n$.
From the identities \eqref{b25}, (a), (b) and \eqref{1prop2} for
$v=J(x,y,z)L_{x,x}^kL_{y,y}^lL_{z,z}^m$ we have:
$$vJ(x,y,z)=$$
$$vG+2/3(-vyzx+vxzy-vzxy+vyxz-vxyz+vzyx-vL_{x,zy}-vL_{z,yx}-vL_{y,xz})=$$
$$J(x,y,z)G L_{x,x}^kL_{y,y}^lL_{z,z}^m-2/3\big(J(x,y,z)yzx+J(x,y,z)xzy-J(x,y,z)zxy+$$
$$J(x,y,z)yxz-J(x,y,z)xyz+J(x,y,z)zyx+J(x,y,z)L_{x,zy}+J(x,y,z)L_{z,yx}+$$
$$+J(x,y,z)L_{y,xz}\big) L_{x,x}^kL_{y,y}^lL_{z,z}^m=J(x,y,z)J(x,y,z)L_{x,x}^kL_{y,y}^lL_{z,z}^m=0$$
That is,  $$J(x,y,z)L_{x,x}^kL_{y,y}^lL_{z,z}^mJ(x,y,z)=0.$$ 
We prove the idetntity $$J(x,y,z)L_{x,y}^nL_{x,x}^kL_{y,y}^lL_{z,z}^mJ(x,y,z)=0.$$
 by induction on $n$. For $n=0$ it is already proven. Suppose that
 $$J(x,y,z)L_{x,y}^iL_{x,x}^kL_{y,y}^lL_{z,z}^mJ(x,y,z)=0,\quad i\leq n.$$
Applying the operator $\Delta_{x}^1(y)$ to this identity and using the induction hypothesis we obtain the required.\\
From this identity by induction on the degree of the operator $L_{y,z}$ using the operator
$\Delta_{z}^1(y)$ we obtain:
$$J(x,y,z)L_{y,z}^pL_{x,y}^nL_{x,x}^kL_{y,y}^lL_{z,z}^mJ(x,y,z)=0.$$
From this identity by induction on the degree of operator $L_{x,z}$ using the operator
$\Delta_{z}^1(x)$ we obtain:
$$J(x,y,z)L_{x,z}^qL_{y,z}^pL_{x,y}^nL_{x,x}^kL_{y,y}^lL_{z,z}^mJ(x,y,z)=0.$$
From this identity by induction on the degree of operator $L_{x,zy}$ using the operator 
$\Delta_{x}^1(zy)$ and the equalities \eqref{2prop2} and  \eqref{3prop2} we obtain:
$$J(x,y,z)L_{x,zy}^sL_{x,z}^qL_{y,z}^pL_{x,y}^nL_{x,x}^kL_{y,y}^lL_{z,z}^mJ(x,y,z)=0,$$
for all $k,l,m,n,p,q,s\in \mathbb{N}\cup \{0\}.$\\

Now we prove \eqref{b27}. First, we prove that
$$J(x,y,z)L_{x,x}^kL_{y,y}^lL_{z,z}^m(J(x,y,z)L_{x,x})=0.$$
From \eqref{a7} and \eqref{2prop2} we obtain:
$$J(x,y,z)L_{x,x}^kL_{y,y}^iL_{y,yx^2}=-J(x,y,z)L_{x,x}^kL_{y,y}^iL_{yx,yx}-J(x,y,z)L_{x,x}^kL_{y,y}^ixL_{y,yx}+$$
$$+J(x,y,z)L_{x,x}^kL_{y,y}^iL_{y,yx}x=-J(x,y,z)L_{x,x}^kL_{y,y}^iL_{yx,yx}+J(x,y,zL_{x,x}^kL_{y,y}^ix)L_{y,yx}=$$
$$-J(x,y,z)L_{x,x}^kL_{y,y}^iL_{yx,yx}.$$
Therefore, the application of the operator $\Delta_{y}^1(yx^2)$ with \eqref{2prop2} and the equality \eqref{b26} give:
$$J(x,y,z)L_{x,x}^kL_{y,y}^lL_{z,z}^m(J(x,y,z)L_{x,x})=0.$$
Induction on the degree of the operator $L_{x,y}$ and application of
$\Delta_{y}^1(x)$ gives:
$$J(x,y,z)L_{x,y}^nL_{x,x}^kL_{y,y}^lL_{z,z}^m(J(x,y,z)L_{x,x})=0.$$
Induction on the degree of the operator $L_{y,z}$ and application of
$\Delta_{z}^1(y)$ gives:
$$J(x,y,z)L_{y,z}^pL_{x,y}^nL_{x,x}^kL_{y,y}^lL_{z,z}^m(J(x,y,z)L_{x,x})=0.$$
Induction on the degree of the operator $L_{x,z}$ and application of
$\Delta_{z}^1(x)$ gives:
$$J(x,y,z)L_{x,z}^qL_{y,z}^pL_{x,y}^nL_{x,x}^kL_{y,y}^lL_{z,z}^m(J(x,y,z)L_{x,x})=0.$$
Induction on the degree of the operator $L_{y,xz}$, application of
$\Delta_{y}^1(xz)$ and  the equality \eqref{b18} give the equality \eqref{b27}.\\

We now prove \eqref{b28}. The following equality is obvious:
$L_{z,yx}\Delta_{z}^1(zx)=L_{zx,yx}$. From  \eqref{3prop2}:
$$J(x,y,z)T_1T_2...T_iL_{zx,zx}=J(x,y,z)T_1T_2...T_iL_{x,x}L_{z,z}-J(x,y,z)T_1T_2...T_iL_{x,z}^2,$$
 for all $T_j\in T, i\in \mathbb{N}\cup \{0\}.$ Applying the operator $\Delta_{z}^1(y)$
to this equality we obtain:
$$J(x,y,z)T_1T_2...T_iL_{zx,yx}=J(x,y,z)T_1T_2...T_iL_{x,x}L_{z,y}-J(x,y,z)T_1T_2...T_iL_{x,z}L_{x,y}.$$
From \eqref{b26} and \eqref{1prop2}  the identity follows:
$$J(x,y,z)T_1T_2...T_mJ(x,y,z)=0,$$ where $T_i\in T\cup \{L_{x,zx}, L_{z,zx} \}.$
The application of the operator $\Delta_{z}^1(zx)$ to this identity, gives:
$$J(x,y,zx)T_1T_2...T_mJ(x,y,z)+J(x,y,z)T_1T_2...T_mJ(x,y,zx)=0,$$ for all $T_i\in T\cup \{L_{x,zx}, L_{z,zx} \}.$
Now write \eqref{b26} as follows:
$$J(x,y,z)L_{z,yx}^sL_{x,z}^qL_{y,z}^pL_{x,y}^nL_{x,x}^kL_{y,y}^lL_{z,z}^mJ(x,y,z)=0,$$
for any $k,l,m,n,p,q,s\in \mathbb{N}\cup \{0\}.$\\
Applying to it the operator $\Delta_{z}^2(zx)$, grouping
corresponding terms, and considering the previous identity, we obtain \eqref{b27}.
 $\blacksquare$\\

{\bf Lemma 6.} In the algebra $M_3$ the following relations hold:\\
\begin{equation}\label{b29}
J(y,z,x(R_{zy}R_x)^n)=J(y,z,x)L_{x,zy}^n,
\end{equation}
\begin{equation}\label{b30}
J(x,y,z)G^n=6^nJ(x,y,z)L_{x,zy}^n.
\end{equation}

{\sc Proof.} Prove first the identity
$$J(x,x(R_{zy}R_x)^nR_{zy},y)=J(x,y,z)yxL_{x,zy}^n$$
by induction on $n$. For $n=0$ this identity is obvious. Suppose that it holds for $n=k$. 
Prove it for $n=k+1$.
 Applying the operator $\Delta_{y}^1(zy)$ to 
the identity $$J(x,txy,y)=J(x,t,y)xy,$$ we obtain
$$J(x,tx(zy),y)+J(x,txy,zy)=J(x,t,zy)xy+J(x,t,y)x(zy).$$ 
If $t=x(R_{zy}R_x)^kR_{zy},$
then $$J(x,tx(zy),y)=-J(x,txy,zy)+J(x,t,y)x(zy).$$ 
And from the identity
\eqref{a9} we obtain:
$$J(x,tx(zy),y)=-J(x,txy,zy)+J(x,t,y)x(zy)=1/2J(x,t,y)[R_x,R_{zy}]+J(x,t,y)x(zy)=$$
$$J(x,t,y)L_{x,zy}=J(x,y,z)yxL_{x,zy}^{k+1}.$$
We now prove  the identity \eqref{b29} by induction on $n$. Let $n=1$.
Apply the operator $\Delta_{y}^1(zy)$ to the identity
$$J(y,z,xyx)=J(y,z,x)L_{x,y}.$$ We have
$J(y,z,x(zy)x)=J(y,z,x)L_{x,zy}$. Suppose the identity holds for
$n=k$. 
We prove  for  $n=k+1$. Substituing in the identity \eqref{a9} $a$ by $y$, $b$ by $z$, $c$ by $x$, applying to it the operator
$\Delta_{y}^1(zy)$  and substituing $t$ by $x(R_{zy}R_x)^k$ we have

$$J(z,y,t(zy)x)-J(z,y,t)L_{x,zy}=$$
$$-J(z,y,tyx)z+J(z,y,t)zL_{x,y}+1/2J(zy,t,x)[R_z,R_y]+1/2J(y,t,x)[R_z,R_{zy}]=$$
$$-1/2J(y,t,x)yz^2-1/2J(y,t,x)zyz-J(z,y,t)L_{x,y}z+$$
$$+J(z,y,t)zL_{x,y}+1/2J(zy,t,x)[R_z,R_y]+1/2J(y,t,x)[R_z,R_{zy}].$$ 
Let $t=x(R_{zy}R_x)^k$. From \eqref{b19}
and the induction hypothesis:
$$J(z,y,t(zy)x)-J(z,y,t)L_{x,zy}=-1/2J(y,t,x)yz^2+1/2J(y,t,x)zyz+1/2J(y,t,x)[R_z,R_{zy}].$$
From the proven identity,  \eqref{2prop2}, and the identity \eqref{a14}:
$$J(z,y,t(zy)x)-J(z,y,t)L_{x,zy}=1/2J(zy,y,x)xL_{x,zy}^kxyz^2-1/2J(zy,y,x)xL_{x,zy}^kxzyz+$$
$$+1/2J(zy,y,x)xL_{x,zy}^kx[R_z,R_{zy}]=1/2J(zy,y,x)L_{x,zy}^kxxyz^2-1/2J(zy,y,x)L_{x,zy}^kxxzyz-$$
$$+1/2J(zy,y,x)L_{x,zy}^kxx[R_z,R_{zy}]=1/2J(zy,y,x(zy)x^3(R_{zy}R_x)^{k-1})yz^2-$$
$$-1/2J(zy,y,x(zy)x^3(R_{zy}R_x)^{k-1})zyz-J(zy,y,x(zy)x^3(R_{zy}R_x)^{k-1})(zy)z=0.$$
We now prove the identity \eqref{b30} by induction on $n$. Let $n=1$. From
\eqref{a3}, \eqref{b24} and \eqref{b18} we obtain:
$$G(J(x,y,z),x,y,z)=2/3J(J(x,y,z),y,z)x+2/3J(x,J(x,y,z),z)y+$$
$$+2/3J(x,y,J(x,y,z))z-2/3J(x,y,z)J(x,y,z)=$$
$$2J(x,y,z)(zy)x+2J(x,y,z)(xz)y+J(x,y,z)(yx)z=$$
$$2J(x,y,z)L_{x,zy}+2J(x,y,z)L_{y,xz}+2J(x,y,z)L_{z,yx}=6J(x,y,z)L_{x,zy}.$$
Suppose the identity holds for $n=k$. We prove it for $n=k+1$. 
From the identities \eqref{a3}, \eqref{b18}, \eqref{b26} and \eqref{b29}:
$$6^kG(J(x,y,z)L_{x,zy}^k,x,y,z)=2/3\cdot 6^k(J(J(x,y,z)L_{x,zy}^k,y,z)x+$$
$$J(x,J(x,y,z)L_{x,zy}^k,z)y+J(x,y,J(x,y,z)L_{x,zy}^k)z-J(x,y,z)L_{x,zy}^kJ(x,y,z))=$$
$$2/3\cdot 6^k(J(x(R_{zy}R_x)^k,y,z),y,z)x+J(x,J(x,y(R_{xz}R_y)^k,z),z)y+J(x,y,J(x,y,z(R_{yx}R_z)^k))z)=$$
$$2\cdot 6^k(J(x(R_{zy}R_x)^k,y,z)(zy)x+J(x,y(R_{xz}R_y)^k,z)(xz)y+J(x,y,z(R_{yx}R_z)^k)(yx)z)=$$
$$2\cdot 6^k(J(x(R_{zy}R_x)^k,y,z)(zy)x+J(x(R_{zy}R_x)^k,y,z)(xz)y+J(x(R_{zy}R_x)^k,y,z)(yx)z)=$$
$$2\cdot 6^k(J(x(R_{zy}R_x)^k,y,z)L_{x,zy}+J(x(R_{zy}R_x)^k,y,z)L_{y,xz}+J(x(R_{zy}R_x)^k,y,z)L_{z,yx})=$$
$$2\cdot 6^k(J(x(R_{zy}R_x)^k,y,z)L_{x,zy}+J(x(R_{zy}R_x)^k,y,z)L_{x,zy}+J(x(R_{zy}R_x)^k,y,z)L_{x,zy})=$$
$$6^{k+1}J(x,y,z)L_{x,zy}^{k+1}.\blacksquare $$\\

{\bf Lemma 7.} The set $$\textbf{U}\cup \textbf{U}x\cup
\textbf{U}y\cup \textbf{U}z\cup \textbf{U}xy\cup \textbf{U}xz\cup
\textbf{U}yz$$ generates a linear space $J(M_3,M_3,M_3)$ over the field $F.$\\

{\sc Proof.} We first prove that if
$$w\in \textbf{U}=\{J(x,y,z) G^k L_{x,x}^l L_{y,y}^m L_{z,z}^n L_{x,y}^pL_{x,z}^q L_{y,z}^r|k,l,m,n,p,q,r \in \mathbb{N}\cup \{0\}\},$$
then
\begin{equation}\label{b31}
wG=6wL_{x,zy}=3w(xyz-zyx)=3w(yzx-xzy)=3w(zxy-yxz),
\end{equation}
\begin{equation}\label{b32}
wxyz=1/6wG+wL_{y,z}x-wL_{x,z}y+wL_{x,y}z.
\end{equation}
From the identity \eqref{b22}:
$$G(wL_{x,x},x,y,z)-G(w,x,y,z)L_{x,x}=-2J(x,y,z)(wx)x-2J(x,y,z)(wL_{x,x})$$
On the other hand, from the identity \eqref{a1}:
$$(wx)(J(x,y,z)x)=wJ(x,y,z)xx+xwJ(x,y,z)x+J(x,y,z)xxw.$$
From \eqref{b26}, \eqref{b27}, \eqref{b28}: $$wxJ(x,y,z)x=-(wx)(J(x,y,z)x)=0.$$ 
Thus,
$$G(wL_{x,x},x,y,z)-G(w,x,y,z)L_{x,x}=0$$ and the operators $\Delta_{x_1}^1(x_2)$, $x_i\in X$ give the identities:
$$G(wL_{x_1,x_2},x,y,z)-G(w,x,y,z)L_{x_1,x_2}=0,$$ for all $x_i\in X.$
Since $w=J(x,y,z)T_1T_2...T_n$, where $$T_i\in
T=\{L_{x,x},L_{y,y},L_{z,z}, L_{x,y}, L_{x,z},
L_{y,z},L_{x,zy}\},$$ and $n\in \mathbb{N}\cup \{0\}$ then the equality \eqref{1prop2} implies that
$w$ can be written
$$w=J(x,y,z)L_{x,zy}^kS_1S_2...S_{n-k},$$ where $S_i\in T\backslash \{L_{x,zy}\}.$
Therefore, using the proven identity and the identity
\eqref{b30} and \eqref{b18} we obtain:
$$G(w,x,y,z)=G(J(x,y,z)L_{x,zy}^kS_1S_2...S_{n-k},x,y,z)=$$
$$G(J(x,y,z)L_{x,zy}^k,x,y,z)S_1S_2...S_{n-k}=$$
$$6J(x,y,z)L_{x,zy}^{k+1}S_1S_2...S_{n-k}=6wL_{x,zy}=6wL_{y,xz}=6wL_{z,yx}.$$

From the identities \eqref{b25}, \eqref{b18}, \eqref{b26} and \eqref{b19}:
$$G(w,x,y,z)=2/3(wyzx-wxzy)+2/3(wzxy-wyxz)+2/3(wxyz-wzyx)+$$
$$+2/3wL_{x,zy}+2/3wL_{y,xz}+2/3wL_{z,yx}-J(x,y,z)w=$$
$$2/3(wyzx-wxzy)+2/3(wzxy-wyxz)+2/3(wxyz-wzyx)+2wL_{z,yx}=$$
$$2/3(wyzx-wxzy)+2/3(wzxy-wyxz)+2/3(wxyz-wzyx)+1/3wG.$$ 
That is,
$$wG=wxyz-wzyx+wyzx-wxzy+wzxy-wyxz=wxyz-wzyx-wyzx+2wL_{y,z}x+$$
$$+wxyz-2wxL_{y,z}-wzyx+2wzL_{y,z}+wxyz-2wL_{x,y}z=3wxyz-3wzyx.$$
The remaining equalities in \eqref{b31} can be proved similarly. Further, from \eqref{b31}
and \eqref{b19}:
$$wxyz=1/3wG+wzyx=1/3wG-wzxy+2wzL_{x,y}=$$
$$1/3wG+wxzy-2wL_{x,z}y+2wzL_{x,y}=1/3wG-wxyz+2wL_{y,z}x-2wL_{x,z}y+2wL_{x,y}z,$$ which implies \eqref{b32}.\\
Let $w\in J(M_3,M_3,M_3)$. Proposition 1 implies that there are
 $\alpha_i$ from $F$ for which
$$w=\sum_{i}\alpha_iJ(x,y,z)x_{i1}x_{i2}...x_{ik_i}, x_{ij}\in X, k_i\in \mathbb{N}\cup \{0\}.$$
Therefore, it suffices to show that the polynomials of the form
$J(x,y,z)x_{1}x_{2}...x_{k}, x_{i}\in X,$ where
$k\in\mathbb{N}\cup\{0\}$ are the linear combinations of the
polynomials from the set $$\textbf{U}\cup \textbf{U}x\cup
\textbf{U}y\cup \textbf{U}z\cup \textbf{U}xy\cup \textbf{U}xz\cup
\textbf{U}yz.$$ A proof by induction on $k$ obviously follows
from the identities \eqref{b19} and \eqref{b32}.$\blacksquare$\\

{\sc Proof of the Theorem.} We prove the speciality of the algebra $M_3.$ That is, we show that
free Malcev algebra on the free generators $X$ is isomorphic to the
subalgebra of the algebra $(Alt[X])^{-}$, generated by the set $X$.  
 Let $f$ be the canonical homomorphism $f:M_3 \longrightarrow (Alt[X])^{-}$.
From Lemma 7 it follows that it suffices to show that the set $f(\textbf{U}\cup \textbf{U}x\cup \textbf{U}y\cup \textbf{U}z\cup
\textbf{U}xy\cup \textbf{U}xz\cup \textbf{U}yz)$ is linearly independent in $Alt[X]$. We prove first
the linear independence of the set $f(\textbf{U})$.
As in ~\cite{Il} we define in $Alt[X]$ the following subsets:
$$ W_0=\{(x,y,z)(L^+_{x,y})^{n_1}(L^+_{y,z})^{n_2}(L^+_{z,x})^{n_3}(L^+_{x,x})^{n_4}(L^+_{y,y})^{n_5}(L^+_{z,z})^{n_6}
(L^+_{x,[y,z]})^{n_7}\mid n_i\in \mathbb{N}\cup \{0\}\},$$
$$W_1=\{[w,a]\mid w\in W_0, a\in X\},$$
$$W_2=\{(w,x,y), (w,y,z), (w,z,x)\mid w\in W_0\},$$
$$W=W_0\cup W_1\cup W_2,$$
$$W'=\{w\circ (x,y,z)\mid w\in W_0\}.$$
We prove that for each $u$ from $U$ there exists 
$\alpha\in F, \alpha\neq 0$, such  that  $f(u)=\alpha w$, where $w\in W$. Using the equality \eqref{1prop2} 
we can write $u$ in the form 
$$u=J(x,y,z)T_1...T_{n}T'_1...T'_{m},$$ where $T_i\in \{L_{x,y}, L_{x,z}, L_{y,z}, L_{x,zy} \}$
$\textrm{and } T_i\in \{L_{x,x}, L_{y,y}, L_{z,z}\}.$ Denoting $J(x,y,z)T_1...T_{n}$ by $u_0$ we  have 
$u=u_0T'_1...T'_{m}.$ 
In first we prove the assertion for $u_0=J(x,y,z)T_1...T_{n}$ using the induction on  $n$
 For  $n=0$ it is obvious. Assume that the assertion is proved for $u=J(x,y,z)T_1...T_{n-1}$ and assume that $T_n=L_{x,y}.$
 It is easy to show that then
 $$f: uL_{x,y}\longmapsto 1/2([[f(u),x],y]+[[f(u),y],x])=$$
 $$1/2(f(u)xy-xf(u)y-y(f(u)x)+y(xf(u))+f(u)yx-yf(u)x-x(f(u)y)+x(yf(u)))=$$
 $$-2f(u)L^+_{x,y}-2f(u)L^+_{y,x}.$$
 
The equalities (8), (9) and (24) of  ~\cite{Il}  can be written in the following form 
$$(a,b,a)^+=0,$$
$$aL^+_{b,c}-aL^+_{c,b}=(b,a,c)^+,$$
$$(x_i,(x,y,z)T^+_1...T^+_n,x_j)^+=0,$$ for  $x_i,x_j \in X$ and  $T^+_s\in \{L^+_{x,y}, L^+_{x,z}, L^+_{y,z}, L^+_{x,zy} \},$
where $(a,b,c)^+$ denotes $aL^+_{b,c}.$
 From these relations and from the  induction hypothesis
  it follows that $f(u)L^+_{x,y}=f(u)L^+_{y,x}$.
 Therefore, $$f: uL_{x,y}\longmapsto -4f(u)L^+_{x,y}.$$ Similarly,
 $$f: uL_{x,zy}\longmapsto 4f(u)L^+_{x,[y,z]}.$$ From \eqref{b30} we obtain
 $$f: uG\longmapsto 24f(u)L^+_{x,[y,z]}.$$
 Finally, direct calculation gives for any $t$ from $M_3$:
 $$f: tL_{x,x}\longmapsto [[f(t),x],x]=x(xf(t))+f(t)xx-xf(t)x-x(f(t)x)=$$
 $$-4(f(t)\circ x\circ x-f(t)\circ (x\circ x))=-4f(t)L^+_{x,x}.$$
  
That is,
 $$f: tL_{x,x}\longmapsto [[f(t),x],x]=-4f(t)L^+_{x,x}.$$
 From the proven identities and from results of ~\cite{Il}  the linear independence of $f(\textbf{U})$ follows.\\
 Note that for any $u\in \textbf{U}$ the following holds $$f: uxy\longmapsto 4f(u)L^+_{x,y}+2(f(u),x,y).$$
It means that $f(\textbf{U})=FW_0$ and the  set $\textbf{U}$ is linearly independent.\\ 
 Now suppose that there exist $u_i \in _F(\textbf{U})$ such  that the following relations are fulfilled. 
$$f(u_0+u_1x+u_2y+u_3z+u_4xy+u_5xz+u_6yz)=0.$$
Denoting $f(u_i)=w_i$ we obtain
$$w_0+[w_1,x]+[w_2,y]+[w_3,z]+4w_4L^+_{x,y}+2(w_4,x,y)+4w_5L^+_{z,x}+2(x,z,w_5)+4w_6L^+_{y,z}+$$
$$+2(y,z,w_6)=(w_0+4w_4L^+_{x,y}+4w_5L^+_{z,x}+4w_6L^+_{y,z})+[w_1,x]+[w_2,y]+[w_3,z]+2(w_4,x,y)+$$
$$+2(x,z,w_5)+2(y,z,w_6).$$
Hence, from the linear independence of the set $W$, proven in ~\cite{Il}, we have:
$w_1=w_2=w_3=w_4=w_5=w_6=0$ and $w_0+4w_4L^+_{x,y}+4w_5L^+_{z,x}+4w_6L^+_{y,z}=0$. Therefore,
$w_0=0$ and from proven we have $u_i=0$, $i\in \{0,1,2,3,4,5,6\}$. \\
Obviously, from Lemma 7 it follows that the set $$\textbf{U}\cup \textbf{U}x\cup
\textbf{U}y\cup \textbf{U}z\cup \textbf{U}xy\cup \textbf{U}xz\cup
\textbf{U}yz$$ is a basis of the space $J(M_3,M_3,M_3)$. $\blacksquare$\\

{\bf  Corollary 5.}  Let $S$ be the free algebra of rank three of the variety generated by a
simple seven-dimensional Malcev algebra. The free Malcev algebra $M_3$ of rank three is a subdirect sum of the
free Lie algebra $L$ of rank three and the free algebra $S$.\\

{\sc Proof.}
Let $K$ be the free algebra of the variety generated
Cayley-Dickson algebra over a field $F$ with a set of free generators $X$. 
It is easy to see that the subalgebra of the Malcev algebra $K^-$,
generated by  $X$ is isomorphic to $S$. From  Theorem 1 of ~\cite{Il} we obtain that
$M_3$ is the subalgebra in $K^-\oplus (Ass[X])^-$. The projection $K^-\oplus (Ass[X])^- \longrightarrow K^-$
induces a homomorphism $g_1: M_3\longrightarrow S$, which obviously is  surjective. In addition, it is similar to show that the
homomorphism $g_2: M_3\longrightarrow L$, induced by the projection $K^-\oplus (Ass[X])^- \longrightarrow (Ass[X])^-$
is surjective.$\blacksquare$\\

{\bf  Corollary 6.} The free Malcev algebra  of rank three is special.\\

{\sc Proof.} It follows from corollary 5.$\blacksquare$\\

In ~\cite{PolSh} particularly it was shown that the free algebra of the variety generated by the Cayley-Dickson algebra is prime.
From this result and  corollary 5 it can be easily obtained\\

{\bf  Corollary 7.} The free Malcev algebra of rank three is semiprime.\\

In conclusion the author expresses his gratitude to I. P. Shestakov both for posing the problems and for his valuable remarks. Also 
the author is very grateful to S.V. Pchelintsev for his inestimable aid of preparing this article.\\

\end{document}